\documentclass[10pt, amsmath, amssymb]{amsart}

\usepackage{amsmath}
\usepackage{amsthm}
\usepackage{amssymb}
\usepackage{graphics}
\usepackage{bm}
\usepackage{natbib}

\newcommand{\be}{\begin{equation}}
\newcommand{\ee}{\end{equation}}
\newcommand{\ba}{\begin{eqnarray}}
\newcommand{\ea}{\end{eqnarray}}
\newcommand{\bi}{\begin{itemize}}
\newcommand{\ei}{\end{itemize}}
\newcommand{\bn}{\begin{enumerate}}
\newcommand{\en}{\end{enumerate}}
\newcommand{\bp}{\begin{proof}}
\newcommand{\ep}{\end{proof}}

\newcommand{\wt}{\ensuremath{\widetilde}}

\newcommand{\mr}{\ensuremath{\mathrm}}

\newcommand{\mc}{\ensuremath{\mathcal}}
\newcommand{\mf}{\ensuremath{\mathfrak}}
\newcommand{\ov}{\ensuremath{\overline}}

\renewcommand{\bm}{\ensuremath{\mathbb }}

\newcommand{\dom}[1]{\ensuremath{\mathrm{Dom} ({#1}) }}
\renewcommand{\dim}[1]{\ensuremath{\mathrm{dim} \left( {#1} \right) }}

\newcommand{\ran}[1]{\ensuremath{\mathrm{Ran} ({#1}) }}
\renewcommand{\ker}[1]{\ensuremath{\mathrm{Ker} ({#1}) }}

\newtheorem{thm}{Theorem}

\newtheorem{lemming}{Lemma}
\newtheorem{prop}{Proposition}

\begin{document}

\bibliographystyle{unsrt}

\title{Characterization of the unbounded bicommutant of $C_0 (N)$ contractions }

\author{R.T.W. Martin}

\address{Department of Mathematics \\ University of California- Berkeley\\
Berkeley, CA, 94720 \\
phone: +1 510 642 6919 \\ fax: +1 510 642 8204}

\email{rtwmartin@gmail.com}

\begin{abstract}
    Recent results have shown that any closed operator $A$
commuting with the backwards shift $S^*$ restricted to $K ^2 _u := H^2 \ominus u H^2$,
where $u$ is an inner function, can be realized as a Nevanlinna function of $S^* _u := S^* | _{K^2 _u}$,
$A = \varphi (S^* _u)$, where $\varphi$ belongs to a certain class of Nevanlinna
functions which depend on $u$. In this paper this result is generalized to show that given any contraction
$T$ of class $C_0 (N)$, that any closed (and not necessarily bounded) operator $A$ commuting with
the commutant of $T$ is equal to $\varphi (T)$ where $\varphi $ belongs to a certain class
of Nevanlinna functions which depend on the minimal inner function $m_T$ of $T$.
\end{abstract}

\maketitle

\section{Introduction}
\label{section:intro}

    Let $u$ be an inner function and let $K^2 _u := H^2  \ominus u H^2$. Recall that the Nevanlinna class $\mc{N}$
in $\bm{D}$ is the class of functions $\varphi = \psi / \chi$ where $\psi , \chi \in H^\infty$
and $\chi$ is not the zero function. The Smirnov class $\mc{N}^+ \subset \mc{N}$ consists of all $\varphi  = \psi /\chi \in \mc{N}$
for which $\chi $ is outer. As defined in \cite{Sarason-ubcom}, the local Smirnov class $\mc{N}^+ _u$ consists
of all $\varphi \in \mc{N}$ for which $u , \chi $ are relatively prime. As discussed in Sects. 3 and 5 of \cite{Sarason-ubcom},
any $\varphi \in \mc{N}^+ _u$ has a unique canonical representation  $\varphi = b/v a$ where $a,b \in H^\infty$, $a$ is an outer
function such that $a(0) = 0$, $|a| ^2 + | b | ^2 =1$ almost everywhere on $\bm{T}$, $v$ is inner and $v,b$ and $v,u$ are relatively
prime. Given $u$ and $K^2 _u$, define the compression $S_u := P_u S | _{K^2 _u}$, where $S$ is the shift (multiplication by $z$) and
$P_u$ is the orthogonal projection of $H^2 $ onto $K^2 _u$. Since $K^2 _u$ is invariant for the backwards shift $S^*$,
$S_u$ is the adjoint of $S^* _u := S^* | _{K^2 _u}$.

Given any $\chi \in H^\infty$ such that $\chi , u$ are relatively prime, one can show that
$\chi  (S _u)$ is injective and has dense range so that $\chi  (S _u) ^{-1}$ can be realized as
a densely defined and closed operator in $K^2 _u$ \cite{Sarason-ubcom} (actually,
the results of \cite{Sarason-ubcom} are expressed in terms of $S^* _u$, we restate them here
in terms of $S_u$).  Hence, as discussed at the end of Sect. 5 of the same paper,
for any $\varphi \in \mc{N} ^+ _u$, one can naturally define $\varphi  (S _u) = \left( (va)  (S _u) \right) ^{-1} b (S _u)$
as a closed operator on a dense domain in $K^2 _u$.

    In \cite{Sarason-ubcom}, Sarason extends the results of Su\'{a}rez in \cite{Suarez} to prove the
following:

\begin{thm}{ (Sarason) }
    A closed operator $A$ densely defined in $K ^2 _u$ commutes with $S _u$ if and only if
$A = \varphi (S _u)$ where $\varphi  \in \mc{N}^+ _u$.  \label{thm:ubcom}
\end{thm}

    This is a natural extension of the following well-known fact, first established in \cite{Sarason-bcom}:

\begin{thm}{ (Sarason) }
    A bounded operator $B$ belongs to the commutant of $S_u$ if and only if $B = h(S_u)$ for some $h
\in H^\infty$.  \label{thm:bcom}
\end{thm}

 Recall that a contraction $T$ is said to be of class $C_0$ if there is
an $H^\infty $ function $v$ such that $v(T) =0$. For any such contraction there is a minimal inner function
$m_T \in H^\infty$ such that $m_T (T)=0$ and $m_T$ is a divisor of any $h \in H^\infty $ for which $h(T) =0$.
The multiplicity $\mu _T $ of $T$ is the minimum cardinal number of a subset $\mf{S} \subset \mc{H}$ such
that $\bigvee  _{n=0} ^\infty T^n \mf{S} = \mc{H}$.

In Sect. 4 of \cite{Nagy-jform}, Sz.-Nagy and Foias use their canonical Jordan model for any contraction $T$
of class $C_0$ with $\mu _T , \mu _{T^*} < \infty$ to show that any element in the double commutant of $T$
is a Nevanlinna function of $T$. Here the double commutant, $(T)^{''}$, is defined, as usual, as the set of all bounded
operators commuting with the commutant, $(T)'$, the set of all bounded operators commuting with $T$.

\begin{thm}{ (Sz.-Nagy, Foias) }
    For any contraction $T$ of type $C_0$ with finite multiplicities $\mu _T <\infty$ and $\mu _{T^*} < \infty$,
all operators $A \in (T) ^{''} $ have the form $\varphi (T)$ where $\varphi \in \mc{N}_T$. \label{thm:bicom}
\end{thm}

They further show by example that there exist such contractions $T$ for which there are $\varphi (T) \in (T) ^{''}$ where
$\varphi \in \mc{N} $, $\varphi \notin H^\infty$ so that $H^\infty $ functions of $T$ do not exhaust
the double commutant of $T$.

Recall that a contraction $T$ is said to be of class $C _0 (N)$ if $T^n $, $(T^*) ^n$ converge strongly to $0$,
and if $N= \mf{d} _T = \mf{d} _{T^*}$, where the deficiency index $\mf{d} _T $ is defined as $\mf{d} _T := \dim{
 \ov{\left( I - T^*T \right) \mc{H}}}$. A contraction $T$ belongs to the class $C_0 (N)$ if and only if it is unitarily equivalent
to some $S (\Theta )$ where $S (\Theta )$ is the compression of the the shift on $H^2 (\mc{H} _N)$, to the subspace
$K^2 (\Theta) := H^2 (\mc{H} _N ) \ominus \Theta H^2 (\mc{H} _N)$. Here $\mc{H} _N$ is any $N$ dimensional Hilbert space,
$H^2 (\mc{H} _N )$ is the Hardy space of functions on the unit disc which take values in the Hilbert space $\mc{H} _N$,
and $\Theta $ is an $N \times N$ matrix valued inner function. The class $C_0 (N)$ is contained in the class of $C_0$ contractions
with finite multiplicities $\mu _T ,\mu _{T^*} < \infty$. Indeed, if $T \in C_0 (N)$ then $\mu _T , \mu _{T^*} \leq N$
(\cite{Nagy-qa}, \cite{Nagy-jform} pg. 94). Moreover, if $T \in C_0 (N)$, and $\Theta _T$ is the $N\times N$ matrix valued inner function such that
$T$ is unitarily equivalent to $S (\Theta _T )$, then $m_T$ is equal to the quotient of $\det (\Theta _T )$ by the greatest
common inner divisor of the minors of order $N-1$ of the matrix of $\Theta _T$ (\cite{Foias}, Chapter VI, Theorem 5.2).

Given a bounded operator $B$, a closed operator $A$ (not necessarily bounded) will be said to commute
with $T$ provided $B: \dom {A} \rightarrow \dom{A}$ and $[A,B] f= (AB -BA) f= 0$ for all $f\in \dom{A}$.
This implies that $AB$ is an extension of $BA$ (in general a proper one) and will be written more concisely as
$AB \supset BA $. Given a contraction $T$, we will say that
a closed operator $A$ belongs to the unbounded double commutant of $T$, $(T) ^{''} _{ub}$,
if $AB \supset BA $ for any $B \in (T) '$. Note that $(T) ^{''} \subset (T) _{ub} ^{''}$.
Just as Sz.-Nagy and Foias used Sarason's original result, Theorem \ref{thm:bcom}
to prove Theorem \ref{thm:bicom}, in this paper we will perform the necessary modifications to the methods
of \cite{Nagy-jform} and use Sarason's new, `unbounded' version, Theorem \ref{thm:ubcom}, of Theorem
\ref{thm:bcom} to prove the following `unbounded' analog of Theorem \ref{thm:bicom}.

\begin{thm}
    Let $T$ be a contraction of class $C_0 (N)$. Then
$A \in (T) ^{''} _{ub}$ if and only if $A = \varphi (T)$ for some $\varphi \in \mc{N}^+ _{m_T}$.
\label{thm:ubicom}
\end{thm}

   Note that our assumptions on $T$ in Theorem \ref{thm:ubicom} are stronger than those used by Sz.-Nagy and
Foias
in Theorem \ref{thm:bicom}.  We expect that there is a stronger version of Theorem \ref{thm:ubicom} which
holds for all $T$ satisfying the conditions of Theorem \ref{thm:bicom}, but this will not be proven here.
The reason for the more restrictive assumption is that our proof will require the use of a lemma that
states that any $T \in C_0 (N)$ cannot be quasi-similar to its restriction to any proper invariant subspace,
see Remark \ref{subsubsection:facts}. If this lemma could be shown to hold for all contractions satisfying the conditions of Theorem \ref{thm:bicom},
\emph{i.e.} all $C_0$ contractions $T$ with finite multiplicities $\mu _T, \mu _{T^*} < \infty$,
then the methods used in this paper would imply that the conclusions of Theorem \ref{thm:ubicom} hold for
this more general class of contractions as well.

\subsection{Contractions of class $C_0$ with finite multiplicity}

    In \cite{Nagy-jform}, a Jordan operator is defined as a contraction of the form:
\be S( u_1, u_2 , ... , u_k ) := S _{u_1} \oplus S _{u_2} \oplus ... \oplus S _{u_k}, \ee where each $u_i$ is a non-constant
inner function and each $u_i$ is an inner divisor of $u_{i-1}$ for $2 \leq i \leq k$. Clearly such
an operator is of class $C_0 (N)$ with minimal function $u_1$.

Recall that a bounded operator $X : \mc{H} _1 \rightarrow \mc{H} _2$ is called a quasi-affinity
if it has dense range, and is injective, \emph{i.e.} if it has a (possibly unbounded) inverse
defined on a dense domain in $\mc{H} _2$. Given $T_i \in B (\mc{H} _i)$, $i=1,2$, $T_1$ is called
a quasi-affine transform of $T_2$ if there exists a quasi-affinity $X$ intertwining $T_2$ and $T_1$,
$T_2 X = X T_1$. This is denoted by $T_2 \succ T_1$.  Note that $\succ$ is a transitive partial order,
and that if $T_1 \succ T_2$ then $T_2 ^* \succ T_1 ^* $. If the $T_i $ are of class $C_0$ and $T_1 \succ T_2$, then
$m_{T_1} = m_{T_2}$, and $\mu _{T_1} \leq \mu _{T_2}$. In particular
if $T_1 \succ T_2$ and $T_2 \succ T_1$, then $T_1, T_2$ are said to be quasi-similar, and $\mu _{T_1} = \mu _{T_2}$.
For more details, we refer the reader to \cite{Nagy-jform}, or \cite{Foias}.

This next theorem of \cite{Nagy-jform} shows that every contraction $T \in C_0$ with finite
multiplicities has a canonical `Jordan model'.

\begin{thm}{ (Sz.-Nagy, Foias) }
    Let $T$ be any contraction of class $C_0$ with finite multiplicities $\mu _T , \mu _{T^*} < \infty$.
Then $T$ is quasi-similar to a Jordan model operator $S(m_1, ..., m_k)$ where $k = \mu _{T} = \mu _{T^*}$,
and $m_1 = m_T$.    \label{thm:jform}
\end{thm}

We will also need the following lemma from \cite{Nagy-jform}.

\begin{lemming}{ (Sz.-Nagy , Foias) }
Let $m, m'$ be non-constant inner functions, with $m$ an inner divisor of $m'$,
$m' = mq$, $q$ inner. Then $\mc{H} ^0 := qH^2 \ominus m' H^2 = K^2 _{m'} \ominus K^2 _q$ is
invariant for $S _{m'}$ and $S^0:= S _{m'} | _{\mc{H} ^0}$ is unitarily equivalent to $S _m$.

Moreover, $R: K^2 _m \rightarrow \mc{H} ^0$, $Rh = qh$ is a unitary transformation
onto $\mc{H} ^0$ such that $S^0 R = RS _m$.

The operator  $Q:= R^{-1} q(S _{m'}) : K^2 _{m'} \rightarrow K^2 _m$, where
$q(S _{m'} ) : K^2 _{m'}  \rightarrow \mc{H} ^0$ is such that $Q K^2 _{m'} = K^2 _m$ and
$S _m Q = Q S _{m'}$. \label{lemming:lem}
\end{lemming}

\section{Proof of Theorem \ref{thm:ubicom}}

    Recall, as defined in \cite{Foias}, for any contraction $T$, $\mc{N}_T$ is the class
of Nevanlinna functions $\varphi = \psi / \chi$ such that $\chi (T)$ is injective and has
dense range. The class of all such $\chi $ is denoted by $K^\infty _T$.
For all $\varphi \in \mc{N}_T$, one can define $\varphi (T) = \chi ^{-1} (T) \psi (T)$.

\subsubsection{Remark} By Theorem 1.1, Chapter 4 of \cite{Foias}, if $\varphi = \psi /\chi  \in \mc{N}_T$, then
$\varphi(T) = \chi ^{-1} (T) \psi (T)$ is a closed operator on a dense domain in $\mc{H}$, and
$\varphi (T) \in (T) ^{''} _{ub}$.

For a contraction of class $C_0$ the following lemma shows that this definition of $\mc{N}_T$
coincides with $\mc{N} ^+ _{m_T}$.

\begin{lemming}
    For a contraction of class $C_0$, $\mc{N}_T = \mc{N} ^+ _{m_T}$. \label{lemming:nevaclass}
\end{lemming}

\begin{proof}
    This lemma is a consequence of known results. If $\varphi \in \mc{N}^+ _{m_T}$, let
$\varphi = b /va $ be its canonical decomposition with $b, v, a \in H^\infty$, $a$, outer,
and $v$ an inner function relatively prime to both $m_T$ and $b$ (see the Introduction, Sect.
\ref{section:intro}, and \cite{Sarason-ubcom}). By Proposition 3.1, Chapter 3, of \cite{Foias},
any outer function belongs to $K^\infty (T)$, the class of $H ^\infty $ functions $\chi$ for which
$\chi (T)$ is injective and has dense range. Secondly, by Proposition 4.7 (b), Chapter 3, of \cite{Foias},
an inner function $v$ belongs to $K^\infty _T$ if and only if $v$ and $m_T$ are relatively prime. Also by
earlier results of the same section (\cite{Foias}, pg. 121), $K^\infty _T$ is multiplicative. It follows
that $va \in K^\infty (T)$ so that $\varphi = b/va \in \mc{N}_T$ and $\mc{N}^+ _{m_T} \subset \mc{N}_T$.

    If $\phi \in \mc{N}_T$, then $\phi $ has the canonical decomposition $\phi = b/ va$ where $v$ is inner and $a$
is outer, and $v,b$ are relatively prime. Using the results mentioned above, and the uniqueness of the
canonical decomposition, it is not difficult to show that if $v, m_T$ are not relatively prime, then
$va  \notin K^\infty _T$, and that $\phi \notin \mc{N}_T$. Hence $v, m_T$ are relatively prime so that
$\mc{N}_T \subset \mc{N}^+ _{m_T}$.
\end{proof}

\begin{prop}
    Let $T$ be a contraction of class $C_0$ with $\mu _T, \mu _{T^*} < \infty$. Then if
$A \in (T) _{ub} ^{''}$, there is a $\varphi \in \mc{N}^+ _T$ and a dense domain $\mf{D} \subset \mc{H}$
such that $A f = \varphi (T) f$ for all $f \in \mf{D}$. \label{prop:neva}
\end{prop}

The following lemma is needed in the proof of Proposition \ref{prop:neva}.

\begin{lemming}{ (Sarason) }
If $A$ is a densely defined operator commuting with every $H^\infty$ function of
$S_u$ then $A$ is closable, and $\ov{A}$ commutes with $S_u$.
\label{lemming:closable}
\end{lemming}

    The above lemma has not been published before. It will appear, along with its
proof, in an upcoming paper by D. Sarason.

The proof of the above proposition follows the proof of Theorem
\ref{thm:bicom} very closely. We will partially sketch the unchanged
portions of the proof, and indicate where the methods of \cite{Nagy-jform} are modified.

\begin{proof}
    By Theorem \ref{thm:jform}, $T$ is quasi-similar to $S:= S(m_1 , ... , m_N) $
acting on $\mc{G} := K^2 _{m_1} \oplus ... \oplus K^2 _{m_N}$ where $N = \mu _T < \infty$
and $m_1 = m_T$. For convenience we will denote $K ^2 _j := K^2 _{m_j}$ and
$S_j := S _{m_j}$. Let $X: \mc{G} \rightarrow \mc{H}$ and $Y: \mc{H} \rightarrow \mc{G}$
be the quasi-affinities such that
\be  TX = XS  \ \ \mr{and} \ \  SY = YT. \label{eq:qsim} \ee

By Lemma \ref{lemming:lem}, for each $j=1, ..., N$ there is a subspace $\mc{H} ^0 _j
\subset K^2 _1$ which is invariant for $S _1$, and such that $S^0 _j :=
S_1 | _{\mc{H} ^0 _j}$ is unitarily equivalent to $S _j$,

\be S^0 _j R_j = R_j S_{m_j} , \label{eq:twine} \ee with each $R_j : K^2 _j \rightarrow \mc{H} ^0 _j$ unitary, $j=1,...,N$. Observe that
$\mc{H} ^0 _1 = K ^2 _{m_1}$. Also by the same lemma there exist bounded operators
$Q_j : K^2 _{m_1} \rightarrow K^2 _{m_j}$ such that
\be S_j Q_j = Q_j S_1 , \label{eq:twiner} \ee where each $Q_j $ is onto $K^2 _{m_j}$.

Now consider the operators
\be S^0 := S^0 _1 \oplus ... \oplus S^0 _N \ \ \mr{on} \ \ \mc{G} ^0 := \mc{H} ^0 _1 \oplus ... \oplus \mc{H} ^0 _N,\ee
and \be \hat{S} := \oplus _{j=1} ^N S_1  \ \ \mr{on} \ \ \hat{\mc{G}} := \oplus _{j=1} ^N K^2 _1.\ee

If we let $R:= \oplus _{j=1} ^N R_j$ and $Q:= \oplus _{j=1} ^N Q_j$ it follows from (\ref{eq:twine})
and (\ref{eq:twiner}) that $R$ is unitary, that
\be S^0 R = RS, \ \ R \mc{G} = \mc{G} ^0, \label{eq:sumtwine} \ee and that
\be SQ = Q \hat{S}, \ \ \  \ Q\hat{\mc{G}} = \mc{G}. \label{eq:sumtwiner} \ee
According to Lemma \ref{lemming:lem}, $Q_j = R _j ^{-1} q _j (S_1 ) $ for the inner function
$q_j$ where $m_1 = m_j q_j$, $1 \leq j \leq N$, and $H^0 _j = K^2 _1 \ominus K^2 _{q_j}$. It
follows that $\ker{Q_j} = K^2 _{q_j}$ so that $K^2 _1 = H^0 _j \oplus \ker{Q_j}$ for each
$1 \leq j \leq N$ and $\hat{\mc{G}} = \mc{G} ^0 \oplus \ker{Q}$.

Using the quasi-affinity $Y$,  (\ref{eq:sumtwine}), and the fact that $\hat{S} | _{\mc{G} ^0} = S^0$
yields
\be \hat{S} R Y = RYT . \label{eq:stot} \ee

    Now choose any $W \in (\hat{S} ) ' $. By (\ref{eq:qsim}), (\ref{eq:sumtwiner}) and
(\ref{eq:stot}), it follows that
\ba T(XQWRY) & = & XSQWRY = XQ\hat{S} WRY = XQW \hat{S} RY \nonumber \\
& = &  (XQWRY) T, \ea so that $XQWRY \in (T) '$. Hence for $A \in (T) ^{''} _{ub}$, one
has $A (XQWRY) \supset (XQWRY) A$. Let $C:= RYXQ$ and define $B:= RYAXQ$ on $\dom{B} := \mf{D} _1 \oplus \ker{Q}
\subset \mc{G} ^0 \oplus \ker{Q} = \hat{\mc{G}}$ where $\mf{D} _1 :=  \{ f \in
\mc{G} ^0 | \ f = RYg ; \ g \in \dom{AXQRY} \}$. The linear manifold $\mf{D} _1$ is clearly dense in $\mc{G} ^0$
 since $I \in (\hat{S} )'$ implies $XQRY \in (T) '$ so that $A (XQRY) \supset (XQRY) A$. This, along
with the facts that $\dom{A}$ is dense in $\mc{H}$  and $RY : \mc{H} \rightarrow \mc{G} ^0$ is a quasi-affinity
imply that $\mf{D} _1 $ is dense in $\mc{G} ^0$, and hence that $\dom{B}$ is dense in $\hat{\mc{G}}$. Observe that the maps $B$, $C$ obey
\be BW C \supset CWB . \label{eq:bcw} \ee
Furthermore, both $B$ and $C$ commute with $\hat{S}$, $[\hat{S} , C ] = 0$ and $B\hat{S} \supset \hat{S} B$. Indeed,
\ba   \hat{S} B & = & \hat{S} (RY) AXQ = (RY)T (A XQ) \subset RY A (TX) Q \nonumber \\
& = & RYA (XS) Q = (RYAXQ) \hat{S} = B\hat{S}  ,\ea the same follows for the bounded operator $C$ when
$A$ is replaced by $I$. Using the same argument as above, and the fact that $A \in (T) ^{''} _{ub}$ further
shows that $h(\hat{S} ) B \subset B h (\hat{S} )$ for any $h \in H^\infty$.

    Since $B, C$ are linear transformations from $\hat{\mc{G}} $ into $\mc{G} ^0 \subset \hat{\mc{G}}$,
they can be viewed as operators on $\hat{\mc{G}}$. Since $W,C$ are bounded we can consider their matrix representations
$W = [ W _{ij} ] $, $C = [C _{ij}]$, $i,j = 1 , ..., N$ with respect to the decomposition $\hat{\mc{G}}= \oplus
_{j=1} ^N K^2 _1$.

We would like to write $B$ as such a matrix. However, since $A$ and hence $B$ is in general unbounded, we need to check that such a matrix representation of $B$ is valid. For example it could be that the domain of $B$ does not contain any vectors of the form
$f = f_1 \oplus 0 ... \oplus 0$ in which case it would not be possible to write $B$ as a matrix with respect to the
decomposition $\hat{\mc{G}} = \oplus _{j=1} ^N K^2 _1$.

 Since $W\in (\hat{S}   )' $, is arbitrary, it can be chosen to be any matrix $W = [W _{ij} ]$, $1 \leq i,j \leq N$,
such that its
entries $W_{ij} \in (S_1) ' $. Now observe that the range of $C= RYXQ$ is dense in $\mc{G} ^0$. This follows
from (\ref{eq:sumtwine}), (\ref{eq:sumtwiner}) and the fact that both $X,Y$ are quasi-affinities.
Choosing $W _j := E_{1j}$, where $E_{ij}$ are the matrix units with respect
to the decomposition $\hat{\mc{G}} = \oplus _{j=1} ^N K^2 _1$, and using $BWC \supset CWB$ (see equation (\ref{eq:bcw})),
it follows that any vector $\oplus _{j=1} ^N \delta _{ij} f_1$ where $\oplus _{j=1} ^N f_j \in C \dom{B}$
belongs to $\dom{B}$ for each $1\leq i \leq N$. Since $\ran{C}$ is dense in $\mc{G}^0 = K^2 _1 \oplus \mc{H} ^0 _2 ... \oplus \mc{H} ^0 _N$,
and $\dom{B}$ is dense in $\hat{\mc{G}}$, it
follows that the set of all such $f_1 $, where $\oplus _{j=1} ^N f_j \in C \dom{B}$
is a dense linear manifold in $K^2 _1$. In conclusion,
\be \dom{B ' } := \{ \oplus _{j=1} ^N f_j \in \dom{B} | \ \oplus _{j=1} ^N \delta _{ij} f_j \in \dom{B} ; \ 1 \leq i \leq N \}
\subset \dom{B} \ee defines a dense linear manifold in $\hat{\mc{G}} = \oplus _{j=1} ^N K^2 _1$. Let $B' := B | _{\dom{B'}}$. It follows that
$B' = [ B' _{ij} ] $ where each $B' _{ij} $ is a densely defined linear operator in $K^2 _1$.
Explicitly, $\dom{B' _{ij}} = \{ f\in K^2 _1 \ |  \ \oplus _{k=1} ^N \delta _{jk} f\in \dom{B'} \}$,
and given $f\in \dom{B' _{ij}}$, $B' _{ij} f= P_i B' f$ where $P_i$ projects $\hat{\mc{G}} =\oplus _{k=1} ^N
K ^2 _1$ onto the i$^{th}$ copy of $K^2 _1$.

We must also check that $WC : \dom{B'} \rightarrow \dom{B'} $ for any $W \in (\hat{S} ) '$
so that we still have \be  CW B ' \subset B' WC, \label{eq:bprimewc} \ee instead of equation (\ref{eq:bcw}).
To do this it suffices to verify that given any
vector $\hat{f} := (f, 0, ..., 0) \in \dom{B'}$ that $WC\hat{f} \in \dom{B'}$. If we choose for example $W =1$,
equation (\ref{eq:bcw}) implies
that $C : \dom{B} \rightarrow \dom{B}$ so that $C\hat{f} = \hat{g} = (g_1, ..., g_N ) \in \dom{B}$. More generally,
$WC \hat{f} = ( \sum _{j=1} ^N W_{1j} g_j , \sum W_{2j} g_j , .... , \sum W_{Nj} g_j )$. To show that
$WC\hat{f} \in \dom{B '}$ we need to show that each $\hat{g} _i = \oplus _{j=1} ^N \delta _{ij} \sum _{k=1} ^N W _{ik} g_k $
belongs to $\dom{B}$, $1 \leq i \leq N$. For example consider the vector $\hat{g} _1 := ( \sum _{j=1} ^N W_{1j} g_j , 0 , ...., 0)$.
This is clearly equal to $W^{(1)} C \hat{f}$ where $W^{(1)} = [ \delta _{1i} W_{ij} ]$ is the matrix obtained by taking
the first row of the matrix representation of $W$ and setting all remaining entries to $0$.
Since $W \in (\hat{S} ) ' $, each $W_{ij}$ commutes with $S_1$,
and so it follows that we also have $W^{(1)} \in (\hat{S} ) ' $. Hence $W^{(1)} C : \dom{B} \rightarrow \dom{B}$,
so that $W ^{(1)} C \hat{f} = \hat{g} _1 \in \dom{B}$. Using similar arguments for the other entries,
it follows that each $\hat{g} _i ; \ 1 \leq i\leq N$ belongs
to $\dom{B}$ so that $ WC\hat{f} \in \dom{B'}$. We conclude that $WC : \dom{B '} \rightarrow \dom{B '} $. This shows that
(\ref{eq:bprimewc}) holds, $B' WC \supset CW B'$. Also it is clear that since $h(\hat{S} ) : \dom{B'} \rightarrow \dom{B'}$,
and the matrix representation of $h(\hat{S}) = h(S_1 ) \oplus ... \oplus h(S_1) $ is diagonal with respect to the decomposition $\hat{G} = \oplus _{j=1} ^N K^2 _1$, we also still have that $B' h (\hat{S} ) \supset h(\hat{S}) B'$ for any $h \in H^\infty$.

Since $C$ commutes with $\hat{S}$ it follows that the matrix entries $C_{ij}$ of $C$ are bounded operators commuting
with $S_1$. By Theorem \ref{thm:bcom}, it follows that there are $H^\infty $ functions $c_{ij}$
such that $C _{ij} = c_{ij} (S_1)$. Similarly, since $B' h( \hat{S} ) \supset h (\hat{S} ) B' $ for
any $h \in H^\infty$ and the matrix representation
of $h (\hat{S}) $ is diagonal, $h(\hat{S}) = h(S_1) \oplus ... \oplus h (S_1) $, it follows that $B' _{ij} h(S_1)
\supset h(S_1) B' _{ij}$. By Lemma \ref{lemming:closable}, each $B' _{ij}$ has a closure $\ov{B' _{ij}}$ commuting
with $S_1$, and by
Theorem \ref{thm:ubcom} and Lemma \ref{lemming:nevaclass}, there exist Nevanlinna functions
$\varphi _{ij} = b _{ij} / \beta _{ij} \in \mc{N}^+ _{m_1} = \mc{N} ^+ _{m_T} = \mc{N}_T$ such that
$\ov{B' _{ij} } = \beta _{ij} ^{-1} (S_1) b_{ij} (S_1)$ where each $\beta _{ij} ^{-1} (S_1 )$ is densely defined
in $K^2 _1$.

Since $B' WC \supset WCB'$ for any $W \in (\hat{S} ) '$, choose $W$ so that $W_{ij} = I$ for $(i,j) = (k,1)$
and all other entries $0$. It follows that there is a dense set $\mc{D} \subset K^2 _1$ such that
\be \beta _{ik} (S_1) ^{-1} b _{ik} (S_1 ) c_{1j} (S_1) f= c _{ik} (S_1) \beta ^{-1} _{1j} (S_1) b_{1j} (S_1) \phi
\ \ \ \ \forall \ f\in \mc{D}. \ee Hence if $d_{ij} := \beta _{ij} c_{ij}$, then $ b_{ik} (S_1 ) d _{1j} (S_1) f= d_{ik} (S_1) b _{1j} (S_1 ) f\ \ \ \forall \ f\in \mc{D} $. Since $d_{ij} , c_{ij} \in H^\infty$
and $\mc{D}$ is dense it follows that
\be b_{ik} (S_1) d_{1j} (S_1 ) = d_{ik} (S_1) b_{1j} (S_1). \label{eq:key} \ee

Since $C$ has dense range in $\mc{G} ^0 = K^2 _1 \oplus _{j=2} ^N \mc{H} ^0 _j$, elements of the form $\sum _{j=1} ^N c_{1j} (S_1) g_j $ where
each $g_j \in K^2 _1$ are dense in $K^2 _1$. Furthermore, since each $\beta _{ij} (S_1)$ has dense range in $K^2 _1$,
it follows that elements of the form $\sum _{j=1} ^N d_{1j} (S_1 ) g_j $, $g_j \in K^2 _1$ are dense in $K^2 _1$.

The proof now proceeds as in \cite{Nagy-jform}, pgs. 109 -111, and leads to the conclusion that there are $v, w \in H^\infty$
such that $v$ and $m_1 = m_T$ are relatively prime and
\be b_{ik} (S_1) v (S_1) - d_{ik} (S_1) w (S_1) = (b_{ik} v) (S_1) - (\beta _{ik} c_{ik} w ) (S_1) =  0 , \ \ 1 \leq i,k \leq N. \ee
Hence, for all $f\in \ran{\beta _{ik}  (S_1) } $ it follows that
\be \left( \beta _{ik} ^{-1} (S_1) b _{ik} (S_1) v(S_1) - c_{ik} (S_1) w(S_1) \right) f=
\left( \ov{B' _{ik} } v(S_1) - C _{ik} w(S_1) \right) f= 0. \ee
In the above note that each $(b_{ik} v) (S_1)$ maps $\dom{ \beta ^{ik} (S_1 ) ^{-1} } \subset \dom {\ov{B' _{ik}}}$ to
itself. It follows that
\ba  0 & = & \left( \ov{B'} v (\hat{S} ) - C w (\hat{S} )  \right) f\nonumber \\
& = & RY \left( AXQ v(\hat{S} ) - XQ w (\hat{S}) \right) f\ \ \ \ \forall \ f\in \dom{B'}.\ea
Since $R$ is unitary, and $Y$ is injective,
\be 0= \left( A v(T) - w(T)  \right) XQ f\ \ \ \ \forall f\in \dom{B'}. \ee
Here, note that $f\in \dom{B'}$, and $B' f= RYAXQ f$ so that $XQ f\in \dom{A}$. Since $A v(T) \supset v(T) A$,
we can now conclude that
\be \left( v(T) A - w(T) \right) g = 0, \ee for all $g \in XQ \dom{B'}$. Since $\dom{B'} $ is dense in $\hat{\mc{G}}$
we have $\ov{XQ \dom{B'}} = \ov{XQ \hat{\mc{G}}} = \ov{X \mc{G}} = \mc{H}$. Hence there is a dense linear manifold of vectors
$\mf{D} := XQ \dom{B'} $ such that for all $f \in \mf{D}$, $ (v(T) A - w(T) ) f =0$, and hence $A f = v(T) ^{-1} w(T) f
= \varphi (T) f $  for all $f \in \mf{D} \subset \dom{\varphi(T)}$ where $\varphi \in \mc{N}_T$ is a Nevanlinna function.
\end{proof}

    The main result, Theorem \ref{thm:ubicom} will now follow from the above proposition once it is established that given any $\varphi \in \mc{N} _T$,
the closed densely defined operator $\varphi (T)$ has no proper closed restrictions or extensions.

\begin{prop}
    Given any contraction $T$ of class $C_0 (N)$, and any $\varphi \in \mc{N}_T$, the closed
operator $\varphi (T)$ has no proper closed densely defined extensions or restrictions belonging to
$(T)^{''} _{ub}$. \label{prop:norest}
\end{prop}

This proposition is a consequence of the Jordan model for $T$ and the following lemma.

\begin{lemming}
    Given an inner function $u$ and $\varphi \in \mc{N} ^+ _u$, the closed operator $\varphi (S _u)$ has
no proper closed densely defined extension or restriction commuting with $S_u$. \label{lemming:norest}
\end{lemming}

    The proof of this lemma follows immediately from the first part of the proof of Lemma 5.7 of \cite{Sarason-ubcom}.

\subsubsection{Remark} The proof of Proposition \ref{prop:norest} relies on the above lemma, as well as
the following two facts taken from \cite{Nagy-jform}. First, if $T , T' $ are contractions of class $C_0$ with finite multiplicities,
then $T \succ T'$ implies $T, T'$ have the same Jordan models. Conversely $T, T'$ having the same Jordan model
implies they are in fact quasi-similar. Furthermore, if $T \in C_0 (N); \ N \geq 1$, then the restriction $T|_S$
to a proper invariant subspace $S \subset \mc{H}$ cannot have the same Jordan model as $T$. That is, $T$ cannot
be quasi-similar to its restriction to a proper invariant subspace. These facts are contained in Corollaries
$1-2$ of \cite{Nagy-jform}. \label{subsubsection:facts}

\begin{proof}{ (Proposition \ref{prop:norest})}

    Let $S_T$ be the Jordan model of $T$ on $\mc{G} := K^2 _1 \oplus ... \oplus K^2 _N$, and as in the proof of
Proposition \ref{prop:neva}, let $X,Y$ be the quasi-affinities such that $X: \mc{G} \rightarrow \mc{H}$, $Y : \mc{H} \rightarrow
\mc{G}$ and $TX = XS_T$, $S_TY = YT$.

Now consider $\varphi (T)$, and suppose that $R$ is a densely defined proper closed restriction of $\varphi (T)$ such that
$R \in (T) ^{''} _{ub}$, and let $\Gamma _\varphi \supsetneqq \Gamma _R$ denote the graphs of these two operators. Note that
$\Gamma _\varphi  , \Gamma _R \subset \mc{H} \oplus \mc{H}$ are invariant for $T\oplus T$, and since
$XY \in (T)'$ it follows that $\Gamma _\varphi , \Gamma _R$ are also invariant for $XY \oplus XY$.

It is straightforward to verify that $\ov{ (Y \oplus Y ) \Gamma _R } $ and $\ov{ (Y\oplus Y ) \Gamma _\varphi}$
are invariant for $S_T \oplus S_T$. I claim that $\ov{(Y\oplus Y) \Gamma _\varphi } \supsetneqq \ov{(Y\oplus Y) \Gamma _R}$.
Suppose to the contrary that $\ov{(Y\oplus Y) \Gamma _\varphi } = \ov{(Y\oplus Y) \Gamma _R}$, so that
$\ov{(XY\oplus XY) \Gamma _\varphi } = \ov{(XY\oplus XY) \Gamma _R} \subset \Gamma _R\subsetneqq \Gamma _\varphi$.
It is clear that both of these subspaces are invariant
for $T \oplus T$. Let $\Pi _1 := T \oplus T | _{\Gamma _\varphi }$ and $\Pi _2 := T \oplus T | _{\ov{(XY\oplus XY) \Gamma _R}}$.

Since $T\oplus T $ is a contraction of class $C_0 (2N)$ on $\mc{H} \oplus \mc{H}$, and $\Gamma _\varphi \subset \mc{H} \oplus
\mc{H}$ is invariant for $T\oplus T$, it follows from Lemma $3.1$, Chapter $IX$ of \cite{Foias} that $\Pi _1 :=
T \oplus T | _{\Gamma _\varphi} $ is a contraction of class $C_0 (N')$, $N' \leq 2N$.

Moreover since by assumption $XY \oplus XY : \Gamma _\varphi \rightarrow \ov{(XY \oplus XY) \Gamma _\varphi }
= \ov{(XY \oplus XY) \Gamma _R}$, it follows that $XY \oplus XY$ is a quasi-affinity such that
$\Pi _2 (XY \oplus XY ) = (XY \oplus XY ) \Pi _1$. This shows
that $\Pi _2 \succ \Pi _1$. By Corollary 1, pg 91 of \cite{Nagy-jform} (see Remark \ref{subsubsection:facts} above),
the $\Pi _i$ are quasi-similar. Since $\Pi _1$ is the restriction of $\Pi _2 $ to the non-trivial
invariant subspace $\Gamma _R \subset \Gamma _\varphi$,
this contradicts Corollary 2, pg 92 of \cite{Foias} (again, see Remark \ref{subsubsection:facts} above),
that no contraction of class $C_0 (N)$ can be
quasi-similar to its restriction to a proper invariant subspace. This contradiction proves that
$\ov{(Y\oplus Y) \Gamma _R} \subsetneqq \ov{(Y\oplus Y) \Gamma _\varphi}$. As remarked earlier,
$\ov{(Y\oplus Y) \Gamma _\varphi} \subset \Gamma (\varphi (S_T) )$ is invariant for $S_T \oplus S_T$.
Hence $\Gamma (R') := \ov{(Y\oplus Y) \Gamma _R} \subsetneqq \Gamma (\varphi (S_T))$ is the graph of a densely defined
closed operator $R'$ which is a
non-trivial proper restriction of $\varphi (S_T)$, and
which commutes with $S_T$. Since $\varphi (S_T) = \varphi (S_1) \oplus ... \oplus \varphi (S_N) $, it follows that
$R' = R'_1 \oplus ... \oplus R' _N$ where each $R' _i $ is a closed restriction of $\varphi (S_i)$ commuting with
$S_i$. Since $R'$ is proper, one of the $R' _i$ must be a proper closed restriction of $\varphi (S_i )$ commuting
with $S_i$. This contradicts Lemma \ref{lemming:norest}, and proves that $\varphi (T)$ has no proper closed restriction
which belongs to $(T) ^{''} _{ub}$.

    Conversely, if $(T) ^{''} _{ub} \ni A \varsupsetneqq \varphi (T)$ is a proper closed extension of $\varphi (T)$,
then $A^* \subsetneqq \varphi (T) ^* = \wt{\varphi} (T ^*)$, where $\wt{\varphi} (z)  := \ov{\varphi (\ov{z})}
\in \mc{N} ^+ _{T^*}$ (see part (v) of Theorem 1.1 in Chapter IV of \cite{Foias}). Since
$T \in C_0 (N)$ implies that $T^* \in C_0 (N)$, the above arguments show that this is not possible.
\end{proof}

    We now have collected all the ingredients needed in the proof of Theorem \ref{thm:ubicom} which we restate below
for convenience.

\setcounter{thm}{3}
\begin{thm}
    Let $T$ be a contraction of class $C_0 (N)$. Then
$A \in (T) ^{''} _{ub}$ if and only if $A = \varphi (T)$ for some $\varphi \in \mc{N} _T = \mc{N}^+ _{m_T}$.
\end{thm}

\begin{proof}
    By Proposition \ref{prop:neva}, there is a $\varphi \in \mc{N}^+ _{m_T}$ and a dense domain
of vectors $\mf{D} \subset \dom{A} \cap \dom{\varphi (T)}$ such that $A f = \varphi (T) f$
for all $f \in \mf{D}$. Let $\dom{A'} := \{ f \in \dom{A} \cap \dom{\varphi (T)} \  |  \
A f = \varphi (T) f \}  \supset \mf{D}$, and define $A' := A | _{\dom{A'}}$. Then $A'$
is a closed restriction of $A$: if $ (f _n ) _{n \in \bm{N} } \subset \dom{A'}$ is such that
$f _n \rightarrow f $ and $A' f _n \rightarrow g$, then $ \varphi (T) f_n
= A f _n \rightarrow g $ so that $ g = \varphi (T) f = A f$ by the fact that
$A, \varphi (T)$ are closed. This proves that $f \in \dom{A'}$ and that $A'$ is closed. Furthermore,
$A' \in (T) ^{''} _{ub} $. To see this, consider arbitrary $W \in (T)'$ and $f \in \dom{A'} \subset \dom{A}$. Since
$W : \dom{A} \rightarrow \dom{A}$, $WA \subset AW$ and $W \varphi (T) \subset \varphi (T) W$, it follows that $ AW f = WA f = WA' f =
W \varphi (T) f = \varphi (T) W f$. This shows that $ A W f = \varphi (T)  Wf$ for all $f \in \dom{A'}$, so that $Wf
\in \dom{A'}$, $WA' f = A' W f$,  $W : \dom {A'} \rightarrow \dom{A'} $ and $W A' \subset A' W $.

Hence $A' \in (T) ^{''} _{ub}$ is a closed restriction of $\varphi (T)$, and Proposition \ref{prop:norest}
implies that $\varphi (T) = A' \subset A$. Since we now have that $A \supset \varphi (T)$, a second application
of Proposition \ref{prop:norest} shows that $A = \varphi (T)$.

\end{proof}


\begin{thebibliography}{1}

\bibitem{Sarason-ubcom}
D.~Sarason.
\newblock Unbounded operators commuting with restricted backwards shifts.
\newblock {\em Oper. Matrices}, 2:583--601, 2008.

\bibitem{Suarez}
D.~Su\'{a}rez.
\newblock Closed commutants of the backwards shift operator.
\newblock {\em Pacific J. Math.}, 179:371--396, 1997.

\bibitem{Sarason-bcom}
D.~Sarason.
\newblock Generalized interpolation in \uppercase{H}$^\infty$.
\newblock {\em Trans. Amer. Math. Soc.}, 127:768--770, 1967.

\bibitem{Nagy-jform}
B.~Sz.-Nagy and C.~Foia\c{s}.
\newblock Mod\`{e}le de \uppercase{J}ordan pour une classe d'op\'{e}rateurs de
  l'espace de \uppercase{H}ilbert.
\newblock {\em Acta Sci. Math.}, 31:91--115, 1970.


\bibitem{Nagy-qa}
B.~Sz.-Nagy and C.~Foia\c{s}.
\newblock Vecteurs cycliques et quasi-affinit\'{e}s.
\newblock {\em Studia Math.}, 31:35--42, 1968.

\bibitem{Foias}
B.~Sz.-Nagy and C.~Foia\c{s}.
\newblock {\em Harmonic analysis of operators on \uppercase{H}ilbert space}.
\newblock American Elsevier publishing company, Inc., New York, N.Y., 1970.

\end{thebibliography}
\end{document}